\documentclass[a4paper,3p]{elsarticle}
\usepackage{amssymb}
\usepackage{amsmath}
\usepackage{mathrsfs}
\usepackage[all,2cell]{xy}
\UseAllTwocells

\newcommand{\cat}[1]{\textup{\bf #1}}
\newcommand{\topsyst}[1]{\ensm{\cat{Af\,Sys}(#1)}}
\newcommand{\topsysto}[1]{\ensm{\cat{Af\,Sys}_0(#1)}}
\newcommand{\loc}[1]{\ensm{\opm{\cat{#1}}}}
\newcommand{\semilatc}[1]{\ensm{\cat{CSLat}(#1)}}

\newcommand{\algo}{\ensm{\cat{Alg}(\Omega)}}
\newcommand{\set}{\ensm{\cat{Set}}}
\newcommand{\topt}[1]{\ensm{\cat{Af\,Spc}(#1)}}
\newcommand{\btop}[1]{\ensm{{#1}\text{-}\cat{Top}}}
\newcommand{\btopsys}[1]{\ensm{{#1}\text{-}\cat{TopSys}}}

\newcommand{\arw}[1]{\ensm{\xrightarrow{#1}}}
\newcommand{\incl}[3]{\ensm{\xymatrix{{#1\,}\ar@{^{(}->}[r]^-{#2} & {#3}}}}
\newcommand{\cell}[4]{\ensm{\xymatrix{#1 \ar@<0.5ex>[r]^-{#2}  \ar@<-0.5ex>[r]_-{#3} &#4}}}

\newcommand{\sqed}{\hfill{\vrule width 3pt height 3pt depth 0pt}}

\newcommand{\ensm}{\ensuremath}
\newcommand{\mcal}[1]{\ensm{\mathcal{#1}}}
\newcommand{\mscr}[1]{\ensm{\mathscr{#1}}}
\newcommand{\mbb}[1]{\ensm{\mathbb{#1}}}
\newcommand{\mfrak}[1]{\ensm{\mathfrak{#1}}}
\newcommand{\opm}[1]{\ensm{{#1}^{op}}}
\newcommand{\seq}{\ensm{\subseteq}}
\newcommand{\unit}{\ensm{\mathfrak{1}}}
\newcommand{\msf}[1]{\ensm{\mathsf{#1}}}
\newcommand{\ovr}[1]{\ensm{\overline{#1}}}
\newcommand{\copwr}[2]{\ensm{{}^{#1}\!#2}}
\newcommand{\indf}[2]{\ensm{\langle{#1}_{#2}\rangle_{n_\lambda}}}
\newcommand{\indfn}[3]{\ensm{\langle{#1}_{#2}\rangle_{#3}}}
\newcommand{\subb}[1]{\ensm{\langle{#1}\rangle}}
\newcommand{\geqs}{\ensm{\geqslant}}

\journal{Fuzzy Sets and Systems}

\begin{document}

\begin{frontmatter}

\title{Sierpinski object for affine systems}

\author[Jeffrey]{Jeffrey T. Denniston}
\address[Jeffrey]{Department of Mathematical Sciences, Kent State University\\
                  Kent, Ohio, USA 44242}
\ead{jdennist@kent.edu}

\author[Austin]{Austin Melton}
\address[Austin]{Departments of Computer Science and Mathematical Sciences, Kent State University\\
                  Kent, Ohio, USA 44242}
\ead{amelton@kent.edu}

\author[Stephen]{Stephen E. Rodabaugh}
\address[Stephen]{College of Science, Technology, Engineering, Mathematics (STEM), Youngstown State University\\
                  Youngstown, Ohio, USA 44555-3347\\}
\ead{serodabaugh@ysu.edu}

\author[Sergey]{Sergey A. Solovyov\tnoteref{GACR}}
\address[Sergey]{Institute of Mathematics, Faculty of Mechanical Engineering, Brno University of Technology, Technick\'{a} 2896/2, 616 69 Brno, Czech Republic}
\ead{solovjovs@fme.vutbr.cz}
\tnotetext[GACR]{This research was supported by the bilateral project ``New Perspectives on Residuated Posets" of the Austrian Science Fund (FWF) (project No. I 1923-N25) and the Czech Science Foundation (GA\v{C}R) (project No. 15-34697L).}

\begin{abstract}
 Motivated by the concept of Sierpinski object for topological systems of S.~Vickers, presented recently by R.~Noor and A.~K.~Srivastava, this paper introduces the Sierpinski object for many-valued topological systems and shows that it has three important properties of the crisp Sierpinski space of general topology.
\end{abstract}


\end{frontmatter}


\newtheorem{thm}{Theorem}
\newtheorem{prop}[thm]{Proposition}
\newtheorem{cor}[thm]{Corollary}
\newtheorem{lem}[thm]{Lemma}
\newdefinition{defn}[thm]{Definition}
\newdefinition{exmp}[thm]{Example}
\newdefinition{rem}[thm]{Remark}
\newdefinition{prob}[thm]{Problem}
\newproof{pf}{Proof}


\section{Introduction}


\par
The notion of \emph{Sierpinski space} $\mcal{S}=(\{0,1\},\{\varnothing,\{1\},\{0,1\}\})$ (see, e.g.,~\cite{Willard2004}) plays a significant role in general topology. In particular, one can show the following three important properties (see, e.g.,~\cite{Adamek2009,Nel1972}):
\begin{enumerate}[(1)]
 \item A topological space is $T_0$ iff it can be embedded into some power of $\mathcal{S}$.
 \item The injective objects in the category of $T_0$ topological spaces are precisely the retracts of powers
       of $\mathcal{S}$.
 \item A topological space is sober iff it can be embedded as a front-closed subspace into some power of $\mathcal{S}$.
\end{enumerate}
Moreover, E.~G.~Manes~\cite{Manes1974,Manes1976} introduced the concept of \emph{Sierpinski object} in categories of structured sets and structure-preserving maps (a subclass of \emph{concrete categories} of~\cite{Adamek2009}), and provided a convenient characterization of the category of topological spaces and continuous maps among such categories in terms of the Sierpinski object (which is precisely the Sierpinski space in the category in question).

\par
Some of the above-mentioned results have already been extended to lattice-valued topology (see, e.g.,~\cite{Lowen1989,Srivastava1984,Srivastava1986}). In particular, there already exists a convenient characterization of the category of fuzzy topological spaces in terms of the Sierpinski object of E.~G.~Manes~\cite{Srivastava1984}.

\par
In~\cite{Vickers1989}, S.~Vickers introduced the concept of \emph{topological system} as a common framework for both point-set and point-free topologies. He showed that the category of topological spaces is isomorphic to a full (regular mono)-coreflective subcategory of the category of topological systems, which gave rise to the so-called \emph{spatialization procedure} for topological systems (from systems to spaces and back). Inspired by the notion of S.~Vickers, R.~Noor and
A.~K.~Srivastava~\cite{Noor2016} have recently presented the concept of Sierpinski object in the category of topological systems, providing topological system analogues of items~(1), (2) above.

\par
Motivated by the notion of lattice-valued topological system of~\cite{Denniston2009,Solovjovs2009} and the results of~\cite{Noor2016}, in this paper, we show lattice-valued system analogues of the above three items (fuzzifying, therefore, some of the achievements of~\cite{Noor2016}). To better incorporate various lattice-valued settings available in the literature, we use the affine context of~Y.~Diers~\cite{Diers1996,Diers1999,Diers2002} and build our systems over an arbitrary variety of algebras (see, e.g.,~\cite{Denniston2016} for the similar approach). Choosing a particular variety gives a particular lattice-valued setting (for example, variety of frames~\cite{Johnstone1982} provides the setting of lattice-valued topological systems of~\cite{Denniston2009}).


\section{Affine spaces and systems}


\par
This section recalls from~\cite{Denniston2016} the notions of affine system and space, and also their related spatialization procedure. To better encompass various many-valued frameworks, we employ a particular instance of the setting of \emph{affine sets} of Y.~Diers~\cite{Diers1996,Diers1999,Diers2002}, which is based in varieties of algebras.

\begin{defn}
 \label{defn:1}
 Let $\Omega=(n_\lambda)_{\lambda\in\Lambda}$ be a family of cardinal numbers, which is indexed by a (possibly, proper or empty) class $\Lambda$. An \emph{$\Omega$-algebra} is a pair $(A,(\omega^A_\lambda)_{\lambda\in\Lambda})$, which comprises a set $A$ and a family of maps $A^{n_\lambda} \arw{\omega^A_\lambda}A$ (\emph{$n_\lambda$-ary primitive operations} on $A$). An \emph{$\Omega$-homomorphism} $(A,(\omega^A_\lambda)_{\lambda\in \Lambda})\arw{\varphi}(B,(\omega^B_\lambda)_{\lambda\in\Lambda})$ is a map
 $A\arw{\varphi}B$, which makes the diagram
 $$
   \xymatrix{A^{n_\lambda} \ar[d]_-{\omega^A_\lambda} \ar[r]^-{\varphi^{n_\lambda}} & B^{n_\lambda}\ar[d]^-{\omega^B_\lambda}\\
             A \ar[r]_-{\varphi} & B\\}
 $$
 commute for every $\lambda\in\Lambda$. \algo\ is the category of $\Omega$-algebras and $\Omega$-homomorp\-hisms, concrete over the category \set\ of sets and maps (with the forgetful functor $|-|$).
 \sqed
\end{defn}

\par
We notice that every concrete category of this paper will use the same notation $|-|$ (which will be not mentioned explicitly) for its respective forgetful functor to the ground category.

\begin{defn}
 \label{defn:2}
 Let \mcal{M} (resp. \mcal{E}) be the class of $\Omega$-homomorphisms with injective (resp. surjective) underlying maps. A \emph{variety of $\Omega$-algebras} is a full subcategory of \algo, which is closed under the formation of products, \mcal{M}-subobjects (subalgebras), and \mcal{E}-quotients (homomorphic images). The objects (resp. morphisms) of a variety are called \emph{algebras} (resp. \emph{homomorphisms}).
 \sqed
\end{defn}

\par
In the following, we provide some examples of varieties, which are relevant to this paper.

\begin{exmp}
 \label{exmp:1}
 \hfill\par
 \begin{enumerate}[(1)]
  \item \semilatc{\bigvee} is the variety of \emph{$\bigvee$-semilattices}, i.e., partially ordered sets, which have
        arbitrary joins.
  \item \cat{Quant} is the variety of \emph{quantales}, i.e., $\bigvee$-semilattices $A$, equipped with a binary operation
        $A\times A\arw{\otimes}A$, which is associative and distributes across $\bigvee$ from both sides, i.e., $a\otimes(\bigvee S)=\bigvee_{s\in S}(a\otimes s)$ and $(\bigvee S)\otimes a=\bigvee_{s\in S}(s\otimes a)$ for every $a\in A$ and every $S\seq A$~\cite{Kruml2008,Rosenthal1990}. \cat{UQuant} is the variety of \emph{unital quantales}, i.e., quantales $A$ having a unit \unit\ for their operation $\otimes$, i.e., $\unit\otimes a=a=a\otimes\unit$ for every $a\in A$.
  \item \cat{Frm} is the variety of \emph{frames}, i.e., unital quantales, for which $\otimes$ is the binary meet operation
        $\wedge$~\cite{Johnstone1982,Picado2012}.
  \item \cat{CBAlg} is the variety of \emph{complete Boolean algebras}, i.e., complete lattices $A$ such that $a\wedge(b\vee
        c)=(a\wedge b)\vee(a\wedge c)$ for every $a$, $b$, $c\in A$, equipped with a unary operation $A\arw{(-)^{\ast}}A$ such that $a\vee a^{\ast}=\top_A$ and $a\wedge a^{\ast}=\bot_A$ for every $a\in A$, where $\top_A$ (resp. $\bot_A$) is the largest (resp. smallest) element of $A$.
        \sqed
 \end{enumerate}
\end{exmp}

\par
We will denote by $\msf{2}$ the two-element algebra $\{\bot_{\msf{2}},\top_{\msf{2}}\}$ in varieties \semilatc{\bigvee}, \cat{Frm}, and \cat{CBAlg}. Also, given an algebra $A$ of a variety \cat{A}, $A\arw{1_A}A$ will stand for the identity map (i.e., homomorphism) on $A$.

\par
From now on, we fix a variety of algebras \cat{A} (for better intuition, one can think of the variety \cat{Frm} of frames, which provides an illustrative example for all the results in this paper).

\begin{defn}
 \label{defn:3}
 Given a functor $\cat{X}\arw{T}\loc{A}$, where \cat{A} is a variety of algebras, \topt{T} is the concrete category over \cat{X}, whose objects (\emph{$T$-affine spaces} or \emph{$T$-spaces}) are pairs $(X,\tau)$, where $X$ is an \cat{X}-object and $\tau$ is an \cat{A}-subalgebra of $TX$; and whose morphisms (\emph{$T$-affine morphisms} or \emph{$T$-morphisms}) $(X_1,\tau_1)\arw{f}(X_2,\tau_2)$ are \cat{X}-morphisms $X_1\arw{f}X_2$ with the property that $\opm{(Tf)}(\alpha)\in\tau_1$ for every $\alpha\in\tau_2$.
 \sqed
\end{defn}

\begin{defn}
 \label{defn:4}
 Given a functor $\cat{X}\arw{T}\loc{A}$, \topsyst{T} is the comma category $(T\downarrow 1_{\loc{B}})$, concrete over the product category $\cat{X}\times\loc{A}$, whose objects (\emph{$T$-affine systems} or \emph{$T$-systems}) are triples $(X,\kappa,A)$, made by \loc{A}-morphisms $TX\arw{\kappa}A$; and whose morphisms (\emph{$T$-affine morphisms} or \emph{$T$-morphisms}) $(X_1,\kappa_1,A_1)\arw{(f,\varphi)}(X_2,\kappa_2,A_2)$ are $\cat{X}\times\loc{A}$-morphisms $(X_1,A_1)\arw{(f,\varphi)}(X_2,A_2)$, which make the following diagram commute
 $$
   \xymatrix{TX_1 \ar[d]_{\kappa_1} \ar[r]^{Tf} & TX_2 \ar[d]^{\kappa_2}\\
             A_1 \ar[r]_{\varphi} & A_2.}
 $$
 \sqed
\end{defn}

\par
In this paper (for the sake of simplicity), we will restrict ourselves to the functor $T$ of the following form.

\begin{prop}
 \label{prop:1}
 Every subcategory \cat{S} of \loc{A} gives rise to a functor $\cat{Set}\times\cat{S}\arw{\mcal{P}_{\cat{S}}}\loc{A}$, $\mcal{P}_{\cat{S}}((X_1,A_1)\arw{(f,\varphi)}(X_2,A_2))=A_1^{X_1}\arw{\mcal{P}_{\cat{S}}(f,\varphi)}A_2^{X_2}$, $\opm{(\mcal{P}_{\cat{S}}(f,\varphi))}(\alpha)=\opm{\varphi}\circ\alpha\circ f$.
\end{prop}

\par
The case $\cat{S}=\{A\arw{1_A}A\}$ provides a functor $\set\arw{\mcal{P}_{A}}\loc{A}$, $\mcal{P}_{A}(X_1\arw{f}X_2)=
A^{X_1}\arw{\mcal{P}_{A}f}A^{X_2}$, $\opm{(\mcal{P}_{A}f)}(\alpha)=\alpha\circ f$. In particular, if $\cat{A}=\cat{CBAlg}$, and $\cat{S}=\{\msf{2}\arw{1_{\cat{2}}}\msf{2}\}$, then one obtains the well-known contravariant powerset functor $\set\arw{\mcal{P}}\loc{CBAlg}$, which is given on a map $X_1\arw{f}X_2$ by $\mcal{P}X_2\arw{\opm{(\mcal{P}f)}}\mcal{P}X_1$ with $\opm{(\mcal{P}f)}(S)=\{x\in X_1\,|\,f(x)\in S\}$. Additionally (following the terminology of~\cite{Rodabaugh1999a}), the case $\cat{S}=\{A\arw{1_A}A\}$ is called \emph{fixed-basis} approach, and all the other instances of \cat{S} are subsumed under \emph{variable-basis} approach.

\par
The following are examples of affine spaces and systems, which are relevant to this paper.

\begin{exmp}
 \label{exmp:2}
 \hfill\par
 \begin{enumerate}[(1)]
  \item If $\cat{A}=\cat{Frm}$, then $\topt{\mcal{P}_{\msf{2}}}$ is the category \cat{Top} of topological spaces.
  \item \topt{\mcal{P}_{A}} is the category $\cat{Af\,Set}(A)$ of affine sets of Y.~Diers. More precisely (as was pointed out by one
        of the referees), the category \topt{\mcal{P}_{A}} is exactly the category $\cat{ASet}(\mathit{\Omega})$ of~\cite{Giuli2009}. Additional studies on some particular cases of the category $\topt{\mcal{P}_{A}}$ can be found in, e.g.,~\cite{Colebunders2010,Giuli2004}. For more details on different variants of affine sets and their respective categories the reader is referred to~\cite{Demirci2014}.
  \item If $\cat{A}=\cat{UQuant}$ or $\cat{A}=\cat{Frm}$, then $\topt{\mcal{P}_{\cat{S}}}$ is the category \btop{\cat{S}} of
        variable-basis lattice-valued topological spaces of S.~E.~Rodabaugh~\cite{Rodabaugh1999a,Rodabaugh2007}.
        \sqed
 \end{enumerate}
\end{exmp}

\begin{exmp}
 \label{exmp:3}
 \hfill\par
 \begin{enumerate}[(1)]
  \item If $\cat{A}=\cat{Frm}$, then \topsyst{\mcal{P}_{\msf{2}}} is the category \cat{TopSys} of topological systems
        of S.~Vickers.
  \item If $\cat{A}=\set$, then \topsyst{\mcal{P}_A} is the category $\cat{Chu}_A$ of Chu spaces over a set $A$ of
        P.-H.~Chu~\cite{Barr1979}.
  \item If $\cat{A}=\cat{Frm}$, then \topsyst{\mcal{P}_{\cat{S}}} is the category \btopsys{\cat{S}} of variable-basis
        lattice-valued topological systems of J.~T.~Denniston, A.~Melton, and S.~E.~Rodabaugh~\cite{Denniston2009,Denniston2012}.
        \sqed
 \end{enumerate}
\end{exmp}

\par
We end this section by providing the promised affine analogue of the system spatialization procedure.

\begin{thm}
 \label{thm:1}
 \incl{\topt{T}}{E}{\topsyst{T},} $E((X_1,\tau_1)\arw{f}(X_2,\tau_2))=(X_1,\opm{e}_{\tau_1},\tau_1)\arw{(f,\varphi)}(X_2,\opm{e}_{\tau_2},\tau_2)$ is a full embedding, where $e_{\tau_i}$ is the inclusion $\tau_i\hookrightarrow TX_i$, and \opm{\varphi} is the restriction $\tau_2\arw{\opm{(Tf)}|_{\tau_2}^{\tau_1}}\tau_1$. $E$ has a right-adjoint-left-inverse $\topsyst{T}\arw{Spat}\topt{T}$,  $Spat((X_1,\kappa_1,B_1)\arw{(f,\varphi)}(X_2,\kappa_2,B_2))=(X_1,\opm{\kappa}_1(B_1))\arw{f}(X_2,\opm{\kappa}_2(B_2))$. \topt{T} is isomorphic to a full (regular mono)-coreflective subcategory of \topsyst{T}.
\end{thm}

\par
We notice that the case $\cat{A}=\cat{Frm}$ and $T=\mcal{P}_{\msf{2}}$ provides the spatialization procedure for topological systems of S.~Vickers, mentioned in the introductory section.

\par
For the sake of convenience, from now on, we will consider the simplest possible case of fixed-basis affine systems. Thus, from now on, we fix an \cat{A}-algebra $L$ (``$L$" is a reminder for ``lattice-valued") and consider the category \topsyst{\mcal{P}_L}, which will be denoted now (for the sake of brevity) \topsyst{L}. Similarly, we will use the notation \topt{L} for the respective category of affine spaces.


\section{Sierpinski object for affine systems}


\par
Motivated by the ideas of R.~Noor and A.~K.~Srivastava~\cite{Noor2016}, in this section, we introduce an affine system analogue of the Sierpinski space. The respective analogue is based in the concept of Sierpinski object in a concrete category of E.~G.~Manes~\cite{Manes1974,Manes1976}. Restated in the modern language of concrete categories of, e.g.,~\cite{Adamek2009}, the concept in question can be defined as follows (cf.~\cite[Definition~2.7]{Noor2016}).

\begin{defn}
 \label{defn:4.1}
 Given a concrete category $\cat{C}$, a \cat{C}-object $S$ is called a \emph{Sierpinski object} provided that for every \cat{C}-object $C$, it follows that the hom-set $\cat{C}(C,S)$ is an initial source (the proof of Proposition~\ref{prop:3} explains the concept of initial source in full detail).
 \sqed
\end{defn}

\par\noindent
This section constructs explicitly Sierpinski object in the category \topsyst{L} (Definition~\ref{defn:5} and Proposition~\ref{prop:3}), which (by analogy with the classical case of topology) is called the Sierpinski ($L$-)affine system.

\par
From now on, we assume that there exists a free \cat{A}-algebra $S$ over a singleton $\mathsf{1}=\{\ast\}$ with the universal map $\mathsf{1}\arw{\eta}|S|$. More precisely, for every \cat{A}-algebra $A$ and every map $\mathsf{1}\arw{f}|A|$, there exists a unique \cat{A}-homomorphism $S\arw{\ovr{f}}A$, making the following triangle commute
$$
  \xymatrix{\mathsf{1} \ar[r]^-{\eta} \ar[rd]_-{f} & |S| \ar@{.>}[d]^-{|\ovr{f}|}\\
                                                   & |A|.}
$$
As can be seen from the discussion at the very end of Section~\ref{sec:1}, all the varieties of Example~\ref{exmp:1} have the required property (for example, if $\cat{A}=\cat{Frm}$, then $S$ is the three-element chain $\{\bot_S,c,\top_S\}$).

\par
We also draw the attention of the reader to the following notational conventions. Given an \cat{A}-algebra $A$ and some element $a\in A$, the unique map $\msf{1}\arw{}|A|$ with value $a$ will be denoted $h_a^A$. Moreover, the product of a set-indexed family of \cat{A}-algebras $\{A_i\mid i\in I\}$ will be denoted $(\prod_{i\in I}A_i\arw{\pi_i}A_i)_{i\in I}$. In particular, if $A_i=A$ for every $i\in I$, then the respective product will be denoted $(A^I\arw{\pi_i}A)_{i\in I}$ (cf., e.g.,~\cite{Adamek2009}).

\begin{defn}
 \label{defn:5}
 \emph{Sierpinski ($L$-)affine system} is the triple $\msf{S}=(|L|,\kappa_S,S)$, in which the map $S\arw{\opm{\kappa}_S}L^{|L|}$ is given by the following diagram
 $$
   \xymatrix{\mathsf{1}\ar[d]_-{h_a^L}\ar[rr]^-{\eta} & & |S| \ar"2,1"_-{|\ovr{h_a^L}|}\ar@{.>}
             [d]^-{|\opm{\kappa}_S|}\\
             |L|&& |L^{|L|}|.\ar[ll]^-{|\pi_a|}}
 $$
 \sqed
\end{defn}

\par
In the case of the category \cat{TopSys}, one gets precisely the Sierpinski object of~\cite{Noor2016}, which has the form $\msf{S}=(|\msf{2}|,\kappa_S,S)$, where $S=\{\bot_S,c,\top_S\}$ and the map $S\arw{\opm{\kappa}_S}\mcal{P}|\msf{2}|$ (recall that $\mcal{P}|\msf{2}|$ stands for the powerset of $|\msf{2}|$) is given by $\opm{\kappa}_S(\bot_S)=\varnothing$, $\opm{\kappa}_S(c)=\{1\}$, and $\opm{\kappa}_S(\top_S)=|\msf{2}|$. For more intuition on Definition~\ref{defn:5}, we consider a quantale-valued example of the Sierpinski affine system, kindly suggested by one of the referees.

\begin{exmp}
 \label{exmp:4}
 Suppose that $\cat{A}=\cat{UQuant}$. We notice first that following~\cite{Rosenthal1990}, the free unital quantale $(S,\otimes,\unit_S)$ over a singleton is the powerset of the set $\mathbb{N}\bigcup\{0\}$ of natural numbers $\{1,2,3,\ldots\}$ and zero, where $\unit_S=\{0\}$, and $\mscr{A}\otimes\mscr{B}=\{n+m\mid n\in\mscr{A},\,m\in\mscr{B}\}$ for every $\mscr{A},\mscr{B}\in S$ (Minkowski addition). In particular, the map $\mathsf{1}\arw{\eta}|S|$ is given by $\eta(\ast)=\{1\}$. Moreover, given a unital quantale $(L,\otimes,\unit_L)$, for every $a\in L$, the map $S\arw{\ovr{h_a^L}}L$ has the form $\ovr{h_a^L}(\mscr{A})=\bigvee\{a^n\mid n\in\mscr{A}\}$, where $a^n=\underbrace{a\otimes\ldots\otimes a}_{\text{$n$-times}}$ for every $n\in\mbb{N}$ and $a^0=\unit_L$. Thus, the map $S\arw{\opm{\kappa}_S}L^{|L|}$ is defined by $(\opm{\kappa}_S(\mscr{A}))(a)=\bigvee\{a^n\mid n\in\mscr{A}\}$.  We recall now that a unital quantale $(L,\otimes,\unit_L)$ is said to be \emph{integral} provided that $\unit_L=\top_L$. Given an integral quantale $L$, for every $a\in L$, it follows that $(a^0=\top_L)\geqs(a^1=a)\geqs(a^2=a\otimes a)\geqs(a^3=a\otimes a \otimes a)\geqs\ldots\geqs(a^n=a\otimes\ldots\otimes a)\geqs\ldots$, since, e.g., $a=a\otimes\top_L\geqs a\otimes a$. As a consequence, $(\opm{\kappa}_S(\mscr{A}))(a)=a^{\bigwedge\mscr{A}}$ for $\mscr{A}\not=\varnothing$, and $(\opm{\kappa}_S(\varnothing))(a)=\bot_L$. Altogether, for an integral quantale $L$, we arrive at the following formula for the map $S\arw{\opm{\kappa}_S}L^{|L|}$, where $1_L$ stands for the identity map on $L$, and, given $a\in L$, $|L|\arw{\underline{a}}L$ denotes the constant map with value $a$:
 $$
   \opm{\kappa}_S(\mscr{A})=
   \begin{cases}
    \underline{\top_L},        & 0\in\mscr{A}\\
    (1_L)^{\bigwedge\mscr{A}}, & 0\not\in\mscr{A}\text{ and }\mscr{A}\not=\varnothing\\
    \underline{\bot_L},        & \mscr{A}=\varnothing.
   \end{cases}
 $$
 Since $0\in\{0,1\}\cap\{0\}$, $\opm{\kappa}_S(\{0,1\})=\underline{\top_L}=\opm{\kappa}_S(\{0\})$, and, therefore, the map $\opm{\kappa}_S$ is not injective.
 \sqed
\end{exmp}

\par
One of the crucial properties (which we will use in the paper) of the Sierpinski system is the following.

\begin{prop}
 \label{prop:2}
 Let $(X,\kappa,A)$ be an affine system.
 \begin{enumerate}[(1)]
  \item For every $a\in A$, there exists a system morphism $(X,\kappa,A)\arw{(f_a,\varphi_a)}\mathsf{S}$ with
        $f_a=\opm{\kappa}(a)$ and $|\opm{\varphi}_a|\circ\eta=h_a^A$.
  \item Every affine morphism $(X,\kappa,A)\arw{(f,\varphi)}\mathsf{S}$ has the form $(f_a,\varphi_a)$ for some $a\in A$.
 \end{enumerate}
\end{prop}
\begin{pf}
 To show item~(1), for every $x\in X$, we consider the following diagram:
 $$
   \xymatrix{\msf{1}\ar@/^1.5pc/"1,4"^-{h_a^A}\ar"3,2"_-{h^L_{(\opm{\kappa}(a))(x)}}\ar[r]^-{\eta} &
             |S|\ar[d]^-{|\opm{\kappa}_S|}
             \ar[rr]^-{|\opm{\varphi}_a|} & & |A| \ar[d]^-{|\opm{\kappa}|}\\
             & |L^{|L|}|\ar[d]^-{|\pi_{f_a(x)}|}\ar[rr]_{|\opm{(\mcal{P}_Lf_a)}|} & & |L^X| \ar[d]^-{|\pi_x|}\\
             & |L| \ar[rr]_-{|1_L|} && |L|.}
 $$
 Commutativity of its outer square implies $|\pi_x\circ\opm{(\mcal{P}_L f_a)}\circ\opm{\kappa}_S|\circ\eta=|\pi_x\circ\opm{\kappa}\circ\opm{\varphi}_a|\circ\eta$, which gives then $\pi_x\circ\opm{(\mcal{P}_L f_a)}\circ\opm{\kappa}_S=\pi_x\circ\opm{\kappa}\circ\opm{\varphi}_a$. As a consequence, $\opm{(\mcal{P}_L f_a)}\circ\opm{\kappa}_S=\opm{\kappa}\circ\opm{\varphi}_a$ as required.

 \par
 To show item~(2), we define $a=|\opm{\varphi}|\circ\eta(\ast)$, which implies $\opm{\varphi}=\opm{\varphi}_a$. Commutativity of the diagram
 $$
   \xymatrix{S\ar[d]_-{\opm{\kappa}_S}\ar[rr]^-{\opm{\varphi}=\opm{\varphi}_a} & & A \ar[d]^-{\opm{\kappa}}\\
             L^{|L|} \ar[rr]_-{\opm{(\mcal{P}_L f)}} && L^X}
 $$
 implies $f(x)\overset{(\dagger)}{=}|\pi_{f(x)}\circ\opm{\kappa}_S|\circ\eta(\ast)=(\opm{\kappa}_S(\eta(\ast)))(f(x))=
 (\opm{(\mcal{P}_Lf)}\circ\opm{\kappa}_S(\eta(\ast)))(x)=(\opm{\kappa}\circ\opm{\varphi}_a(\eta(\ast)))(x)=(\opm{\kappa}(a))(x)$, where $(\dagger)$ relies on the definition of the map $\opm{\kappa}_S$ from Definition~\ref{defn:5}.
 \qed
\end{pf}

\par
We end the section with the second important (for this paper) property of the Sierpinski affine system.

\begin{prop}
 \label{prop:3}
 \msf{S} is a Sierpinski object in \topsyst{L}.
\end{prop}
\begin{pf}
 In view of Definition~\ref{defn:4.1}, we have to show that for every affine system $(X,\kappa,A)$, the source $\mcal{F}=\topsyst{L}((X,\kappa,A),\msf{S})$ is initial. We check thus that every $\set\times\loc{A}$-morphism $|(\tilde{X},\tilde{\kappa},\tilde{A})|\arw{(g,\psi)}|(X,\kappa,A)|$, such that $|(\tilde{X},\tilde{\kappa},\tilde{A})|\arw{(g,\psi)}|(X,\kappa,A)|\arw{|(f,\varphi)|}|\msf{S}|$ is in \topsyst{L} for every $(f,\varphi)\in\mcal{F}$, is itself in \topsyst{L}.

 \par
 By Proposition~\ref{prop:2}, it follows that for every $a\in A$, $(f_a,\varphi_a)\in\mcal{F}$, and, therefore, $(f_a,\varphi_a)\circ(g,\psi)=(f_a\circ g,\varphi_a\circ\psi)\in\topsyst{L}((\tilde{X},\tilde{\kappa},\tilde{A}),\msf{S})$.
  Fix now $a\in A$ and consider the diagram
 $$
   \xymatrix{S \ar[ddd]_-{\opm{\kappa}_S}\ar[rd]^-{\opm{\varphi}_a}\ar[rr]^-{\opm{\psi}\circ\opm{\varphi}_a} &&
             \tilde{A}\ar[ddd]^-{\opm{\tilde{\kappa}}}\\
             &A \ar[d]^-{\opm{\kappa}}\ar[ru]^-{\opm{\psi}}&\\
             &L^X \ar[rd]|{\opm{(\mcal{P}_Lg)}}&\\
             L^{|S|}\ar[ru]|{\opm{(\mcal{P}_Lf_a)}}\ar[rr]_-{\opm{(\mcal{P}_L(f_a\circ g))}}&& L^{\tilde{X}}.}
 $$
 Commutativity of its outer square provides then $\opm{(\mcal{P}_Lg)}\circ\opm{\kappa}\circ\opm{\varphi}_a=
 \opm{\tilde{\kappa}}\circ\opm{\psi}\circ\opm{\varphi}_a$, which implies $|\opm{(\mcal{P}_Lg)}
 \circ\opm{\kappa}\circ\opm{\varphi}_a|\circ\eta=|\opm{\tilde{\kappa}}\circ\opm{\psi}\circ\opm{\varphi}_a|
 \circ\eta$ and, therefore, $|\opm{(\mcal{P}_Lg)}\circ\opm{\kappa}|\circ h_a^A=|\opm{\tilde{\kappa}}\circ\opm{\psi}|\circ h_a^A$. As a consequence, $\opm{(\mcal{P}_Lg)}\circ\opm{\kappa}(a)=\opm{\tilde{\kappa}}\circ\opm{\psi}(a)$ as required.
 \qed
\end{pf}


\section{Properties of the Sierpinski affine system}


\par
In this section, we are going to show affine system analogues of the three properties of the Sierpinski space, mentioned in the introductory section.


\subsection{$T_0$ affine systems}


\par
We begin with an affine modification of the concept of $T_0$ topological system of~\cite{Vickers1989}.

\begin{defn}
 \label{defn:6}
 An affine system $(X,\kappa,A)$ is $T_0$ provided that for every $x,y\in X$, $(\opm{\kappa}(a))(x)=(\opm{\kappa}(a))(y)$ for every $a\in A$ implies $x=y$.
 \sqed
\end{defn}

\par
The following result shows an (possibly, expected by the reader) example of $T_0$ affine systems.

\begin{prop}
 \label{prop:4}
 The Sierpinski affine system $\mathsf{S}$ is $T_0$.
\end{prop}
\begin{pf}
 Given $a,b\in L$ such that $(\opm{\kappa}(c))(a)=(\opm{\kappa}(c))(b)$ for every $c\in S$, it follows that $(\opm{\kappa}(\eta(\ast)))(a)=(\opm{\kappa}(\eta(\ast)))(b)$, which implies $|\pi_a\circ\opm{\kappa}|\circ\eta(\ast)=
 |\pi_b\circ\opm{\kappa}|\circ\eta(\ast)$ and then $h_a^L=h_b^L$, namely, $a=b$.
 \qed
\end{pf}

\par
We need now a couple preliminary results to arrive at the main theorem of this subsection, which will characterize $T_0$ affine systems in terms of the Sierpinski affine system.

\begin{prop}
 \label{prop:5}
 For every $T_0$ affine system $(X,\kappa,A)$, $\mcal{F}=\topsyst{L}((X,\kappa,A),\msf{S})$ is a mono-source.
\end{prop}
\begin{pf}
 We have to show that given two affine system morphisms \cell{(\tilde{X},\tilde{\kappa},\tilde{A})}{(f_1,\varphi_1)}{(f_2,\varphi_2)}{(X,\kappa,A),} if $(f,\varphi)\circ(f_1,\varphi_1)=(f,\varphi)\circ(f_2,\varphi_2)$ for every $(f,\varphi)\in\mcal{F}$, then $(f_1,\varphi_1)=(f_2,\varphi_2)$.

 \par
 To show that $\varphi_1=\varphi_2$, we fix $a\in A$ and notice that $(f_a,\varphi_a)\in\mcal{F}$ (by Proposition~\ref{prop:2}). Thus, $(f_a,\varphi_a)\circ(f_1,\varphi_1)=(f_a,\varphi_a)\circ(f_2,\varphi_2)$, which implies $\varphi_a\circ\varphi_1=\varphi_a\circ\varphi_2$, i.e., $\opm{\varphi}_1\circ\opm{\varphi}_a=\opm{\varphi}_2\circ\opm{\varphi}_a$. Then $|\opm{\varphi}_1\circ\opm{\varphi}_a|\circ\eta=|\opm{\varphi}_2\circ\opm{\varphi}_a|\circ\eta$ gives $\opm{\varphi}_1\circ h^A_a=\opm{\varphi}_2\circ h^A_a$ and then $\opm{\varphi}_1(a)=\opm{\varphi}_2(a)$.

 \par
 To show that $f_1=f_2$, we fix $x\in\tilde{X}$. For every $a\in A$, $(f_a,\varphi_a)\in\mcal{F}$ and thus (by the discussion in the previous paragraph) $f_a\circ f_1=f_a\circ f_2$. It follows then that for every $a\in A$, $(\opm{\kappa}(a))(f_1(x))=f_a(f_1(x))=f_a(f_2(x))=(\opm{\kappa}(a))(f_2(x))$. Since $(X,\kappa,A)$ is $T_0$, one gets that $f_1(x)=f_2(x)$.
 \qed
\end{pf}

\par
As a convenient condition to continue with, from now on, we assume that our fixed variety \cat{A} has coproducts. The reader could notice that all the varieties of Example~\ref{exmp:1} except \cat{CBAlg}~\cite{Adamek2009} have this property. The coproduct of a set-indexed family of \cat{A}-algebras $\{A_i\mid i\in I\}$ will be denoted $(A_i\arw{\mu_i}\coprod_{i\in I}A_i)_{i\in I}$. In particular, if $A_i=A$ for every $i\in I$, then the respective coproduct will be denoted $(A\arw{\mu_i}\copwr{I}{A})_{i\in I}$ (cf.~\cite{Adamek2009}).

\begin{prop}
 \label{prop:6}
 The category \topsyst{L} has products.
\end{prop}
\begin{pf}
 Let $\{(X_i,\kappa_i,A_i)\mid i\in I\}$ be a set-indexed family of affine systems. We show that its product is given by $((\prod_{i\in I}X_i,\kappa,\coprod_{i\in I}A_i)\arw{(\pi_i,\opm{\mu}_i)}(X_i,\kappa_i,A_i))_{i\in I}$, where $\kappa$ is the unique map, defined by the diagram
 $$
   \xymatrix{A_i\ar[d]_-{\opm{\kappa}_i}\ar[rr]^-{\mu_i} && \coprod_{i\in I}A_i \ar@{.>}[d]^-{\opm{\kappa}}\\
             L^{X_i} \ar[rr]_-{\opm{(\mcal{P}_L\pi_i)}} && L^{\prod_{i\in I}X_i}.}
 $$

 \par
 Given a source $((\tilde{X},\tilde{\kappa},\tilde{A})\arw{(f_i,\varphi_i)}(X_i,\kappa_i,A_i))_{i\in I}$ in \topsyst{L}, we have the following two maps (defined by the universal property of (co)products)
 $$
  \begin{matrix}
   \xymatrix{\tilde{X} \ar@{.>}[d]_-{f} \ar[rd]^-{f_i}&\\
             \prod_{i\in I}X_i \ar[r]_-{\pi_i} & X_i}
   &
   \xymatrix{A_i\ar[r]^-{\mu_i} \ar[rd]_-{\opm{\varphi}_i} & \coprod_{i\in I}A_i \ar@{.>}[d]^-{\opm{\varphi}}\\
             &\tilde{A},}
  \end{matrix}
 $$
 i.e., there exists a unique $\set\times\loc{A}$-morphism $|(\tilde{X},\tilde{\kappa},\tilde{A})|\arw{(f,\varphi)}|(\prod_{i\in I}X_i,\kappa,\coprod_{i\in I}A_i)|$, making the diagram
 $$
   \xymatrix{|(\tilde{X},\tilde{\kappa},\tilde{A})| \ar@{.>}[d]_-{(f,\varphi)} \ar"2,3"^-{|(f_i,\varphi_i)|} &&\\
             |(\prod_{i\in I}X_i,\kappa,\coprod_{i\in I}A_i)| \ar[rr]_-{|(\pi_i,\opm{\mu}_i)|} & & |(X_i,\kappa_i,A_i)|}
 $$
 commute. To show that the morphism in question lies in \topsyst{L}, we consider the following diagram
 $$
   \xymatrix{A_i\ar@/^1.5pc/"1,5"^-{\opm{\varphi}_i} \ar[d]_-{\opm{\kappa}_i}\ar[rr]_-{\mu_i} &&  \coprod_{i\in I}A_i
             \ar[d]^-{\opm{\kappa}} \ar[rr]_-{\opm{\varphi}}&& \tilde{A} \ar[d]^-{\opm{\tilde{\kappa}}}\\
             L^{X_i} \ar@/_1.5pc/"2,5"_-{\opm{(\mcal{P}_Lf_i)}} \ar[rr]^-{\opm{(\mcal{P}_L\pi_i)}} && L^{\prod_{i\in I}X_i} \ar[rr]^-{\opm{(\mcal{P}_L f)}} && L^{\tilde{X}}.}
 $$
 Commutativity of its outer square provides then $\opm{(\mcal{P}_Lf)}\circ\opm{\kappa}\circ\mu_i=
 \opm{\tilde{\kappa}}\circ\opm{\varphi}\circ\mu_i$ for every $i\in I$ and, therefore, $\opm{(\mcal{P}_Lf)}\circ\opm{\kappa}=\opm{\tilde{\kappa}}\circ\opm{\varphi}$, which concludes the proof.
 \qed
\end{pf}

\par
As a particular case of Proposition~\ref{prop:6} (for $\cat{A}=\cat{Frm}$ and $L=\msf{2}$), one gets the construction of products of topological systems from~\cite{Vickers1989} (also employed in~\cite{Noor2016}).

\begin{prop}
 \label{prop:7}
 An affine morphism is monic iff its underlying ground category morphism is monic.
\end{prop}
\begin{pf}
 It will be enough to show that given an \topsyst{L}-monomorphism $(X_1,\kappa_1,A_1)\arw{(f,\varphi)}(X_2,\kappa_2,A_2)$, $X_1\arw{f}X_2$ is injective and $A_2\arw{\opm{\varphi}}A_1$ is an \cat{A}-epimorphism.

 \par
 For the first statement, take two maps \cell{Y}{g}{h}{X_1} with $f\circ g=f\circ h$. The lower part of the diagram
 $$
   \xymatrix{L^{X^2}\ar[d]_-{\opm{(\mcal{P}_Lf)}}& \ar[l]_-{\opm{\kappa}_2}A_2\ar[d]^-{\opm{\varphi}}&\\
             L^{X_1} \ar[d]_-{\opm{(\mcal{P}_Lg)}}& \ar[l]_-{\opm{\kappa}_1} A_1 \ar@{.>}[d]^-{\varphi_g}\ar[rd]^-{1_{A_1}}&\\
             L^Y & \ar[l]^-{\pi_{L^Y}}L^Y\times A_1 \ar[r]_-{\pi_{A_1}}& A_1}
 $$
 provides then (by the universal property of products) an \topsyst{L}-morphism $(Y,\opm{\pi}_{L^Y},L^Y\times A_1)\arw{(g,\opm{\varphi}_g)}(X_1,\kappa_1,A_1)$. In a similar way, one obtains an \topsyst{L}-morphism $(Y,\opm{\pi}_{L^Y},L^Y\times A_1)\arw{(h,\opm{\varphi}_h)}(X_1,\kappa_1,A_1)$. To show $\varphi_g\circ\opm{\varphi}=\varphi_h\circ\opm{\varphi}$, we notice that, first,  $\pi_{A_1}\circ\varphi_g\circ\opm{\varphi}=\opm{\varphi}=\pi_{A_1}\circ\varphi_h\circ\opm{\varphi}$ and, second, $\pi_{L^Y}\circ\varphi_g\circ\opm{\varphi}=\opm{(\mcal{P}_Lg)}\circ\opm{(\mcal{P}_Lf)}\circ\opm{\kappa}_2=\opm{(\mcal{P}_L(f\circ g))}\circ\opm{\kappa}_2=\opm{(\mcal{P}_L(f\circ h))}\circ\opm{\kappa}_2=\pi_{L^Y}\circ\varphi_h\circ\opm{\varphi}$, which  implies the desired equality. Thus, $(f,\varphi)\circ(g,\varphi_g)=(f,\varphi)\circ(h,\varphi_h)$, which gives $(g,\varphi_g)=(h,\varphi_h)$, i.e., $g=h$.

 \par
 For the second statement, take two \cat{A}-morphisms \cell{A_1}{\psi}{\theta}{B} with $\psi\circ\opm{\varphi}=\theta\circ\opm{\varphi}$. The right-hand side of the diagram (where $\varnothing\arw{!_{X_1}}X_1$ denotes the unique possible map, and $L^{\varnothing}$ is a terminal object in \cat{A})
 $$
   \xymatrix{A_2\ar[d]_-{\opm{\kappa}_2}\ar[rr]^-{\opm{\varphi}} && A_1 \ar[d]_-{\opm{\kappa}_1}\ar[rr]^-{\psi} &&
             B\ar[d]^-{!}\\
             L^{X_2}\ar[rr]_-{\opm{(\mcal{P}_Lf)}} && L^{X_1}\ar[rr]_-{\opm{(\mcal{P}_L!_{X_1})}} && L^{\varnothing}}
 $$
 provides then an \topsyst{L}-morphism $(\varnothing,!,B)\arw{(!_{X_1},\opm{\psi})}(X_1,\kappa_1,A_1)$. Similarly, one obtains an \topsyst{L}-morphism $(\varnothing,!,B)\arw{(!_{X_1},\opm{\theta})}(X_1,\kappa_1,A_1)$. Since $(f,\varphi)\circ(!_{X_1},\opm{\psi})=(f,\varphi)\circ(!_{X_1},\opm{\theta})$, we arrive at the conclusion that $(!_{X_1},\opm{\psi})=(!_{X_1},\opm{\theta})$, namely, $\psi=\theta$.
 \qed
\end{pf}

\par
Combining the above results, we arrive at the next proposition (we recall from~\cite{Adamek2009} that an \emph{embedding} in a concrete category is an initial morphism, whose underlying ground category morphism is monic).

\begin{prop}
 \label{prop:8}
 Every $T_0$ affine system can be embedded into some power of \msf{S}.
\end{prop}
\begin{pf}
 Given a $T_0$ affine system $(X,\kappa,A)$, by Propositions~\ref{prop:3}, \ref{prop:5}, the source $\mcal{F}=\topsyst{L}((X,\kappa,A),\msf{S})$ is an initial mono-source. By Proposition~10.26\,(1) of~\cite{Adamek2009}, Proposition~\ref{prop:7} of this paper, and easy calculations with initial sources (given a concrete category \cat{C} with products, if a set-indexed source $(C_1\arw{f_i}C_2)_{i\in I}$ is initial, then the unique morphism $C_1\arw{f}C_2^I$, defined by the universal property of products, is initial; the easy proof relies on the following commutative (for every $i\in I$) diagram
 $$
   \xymatrix{|C|\ar[rd]_-{|g|} \ar[r]^-{h} &|C_1| \ar[d]_-{|f|} \ar[rd]^-{|f_i|}&\\
             &|C_2^I| \ar[r]_-{|\pi_i|} & |C_2|,}
 $$
 which implies that $|C|\arw{h}|C_1|$ is a \cat{C}-morphism), it follows that the unique morphism $(X,\kappa,A)\arw{(f,\varphi)}\msf{S}^{\mcal{F}}$, defined by the universal property of products, is an embedding.
 \qed
\end{pf}

\par
To show the opposite direction of Proposition~\ref{prop:8}, we have to do a bit more.

\begin{prop}
 \label{prop:9}
 The product of a set-indexed family of $T_0$ affine systems is $T_0$.
\end{prop}
\begin{pf}
 Given a set-indexed family $\{(X_i,\kappa_i,A_i)\mid i\in I\}$ of $T_0$ affine systems, we show that the product $(\prod_{i\in I}X_i,\kappa,\coprod_{i\in I}A_i)$ is $T_0$. Take $x,y\in\prod_{i\in I}X_i$ such that $(\opm{\kappa}(a))(x)=(\opm{\kappa}(a))(y)$ for every $a\in\coprod_{i\in I}A_i$. Given $i\in I$, for every $b\in A_i$, $(\opm{\kappa}_i(b))(x_i)=(\opm{\kappa}_i(b))\circ\pi_i(x)=(\opm{(\mcal{P}_L\pi_i)}\circ\opm{\kappa}_i(b))(x)=
 ((\opm{\kappa}\circ\mu_i)(b))(x)=((\opm{\kappa}\circ\mu_i)(b))(y)=(\opm{\kappa}_i(b))(y_i)$ and thus $x_i=y_i$ (since $(X_i,\kappa_i,A_i)$ is $T_0$). It follows then that $x=y$.
 \qed
\end{pf}

\begin{prop}
 \label{prop:10}
 Subobjects of $T_0$ affine systems are $T_0$.
\end{prop}
\begin{pf}
 Given an \topsyst{L}-monomorphism $(X_1,\kappa_1,A_1)\arw{(f,\varphi)}(X_2,\kappa_2,A_2)$ such that $(X_2,\kappa_2,A_2)$ is $T_0$, we have to show that $(X_1,\kappa_1,A_1)$ is also $T_0$. Consider the diagram
 $$
   \xymatrix{A_2 \ar[d]_-{\opm{\kappa}_2}\ar[rr]^-{\opm{\varphi}} && A_1\ar[d]^-{\opm{\kappa}_1}\\
             L^{X_2}\ar[rr]_-{\opm{(\mcal{P}_Lf)}} && L^{X_1}}
 $$
 and take $x,y\in X_1$ such that $(\opm{\kappa}_1(a))(x)=(\opm{\kappa}_1(a))(y)$ for every $a\in A_1$. Then $f(x),f(y)\in X_2$ and, for every $b\in A_2$, $(\opm{\kappa}_2(b))(f(x))=(\opm{(\mcal{P}_Lf)}\circ\opm{\kappa}_2(b))(x)=(\opm{\kappa}_1\circ\opm{\varphi}(b))(x)=
 (\opm{\kappa}_1\circ\opm{\varphi}(b))(y)=(\opm{\kappa}_2(b))(f(y))$, which implies $f(x)=f(y)$, and then $x=y$ by Proposition~\ref{prop:7}.
 \qed
\end{pf}

\par
We are now ready to state the main result of this subsection.

\begin{thm}
 \label{thm:2}
 An affine system $(X,\kappa,A)$ is $T_0$ iff it can be embedded into some power of \msf{S}.
\end{thm}
\begin{pf}
 Proposition~\ref{prop:8} provides the necessity. The sufficiency is given by Propositions~\ref{prop:4}, \ref{prop:9}, \ref{prop:10}.
 \qed
\end{pf}


\subsection{\label{subsec:1}Injective $T_0$ affine systems}


\par
This subsection characterizes injective objects (w.r.t. a certain class of morphisms) in the full subcategory \topsysto{L} of \topsyst{L} of $T_0$ affine systems. We begin with some preliminary assumptions and results.

\par
To prove the next proposition (which restores the convenient setting of the category \topsyst{L}), from now on, we assume that, first, every algebra of our variety \cat{A} is non-empty (which obviously holds for all the varieties of Example~\ref{exmp:1}), and, second, our fixed algebra $L$ has more than one element (which does nothing, apart from excluding the trivial case of a singleton).

\begin{prop}
 \label{prop:11}
 In the category \topsysto{L}, an affine morphism is monic iff its underlying ground category morphism is monic.
\end{prop}
\begin{pf}
 In view of Proposition~\ref{prop:7} and its proof, it will be enough to show that both $(Y,\opm{\pi}_{L^Y},L^Y\times A_1)$ and $(\varnothing,!,B)$ are $T_0$ affine systems. While the latter statement is clear, the former requires a small effort. Take $x,y\in Y$ such that $(\pi_{L^Y}(\alpha,a))(x)=(\pi_{L^Y}(\alpha,a))(y)$ for every $(\alpha,a)\in L^Y\times A_1$. Since $A_1$ is non-empty, we obtain $\alpha(x)=\alpha(y)$ for every $\alpha\in L^Y$. Our assumption on $L$ implies then that $x=y$.
 \qed
\end{pf}

\begin{prop}
 \label{prop:12}
 In the category \topsysto{L}, embeddings are exactly the monomorphisms.
\end{prop}
\begin{pf}
 Every embedding is a monomorphism by its very definition. To show the opposite, we fix an \topsysto{L}-monomorphism $(X_1,\kappa_1,A_1)\arw{(f,\varphi)}(X_2,\kappa_2,A_2)$. By Proposition~\ref{prop:11}, we get that $A_2\arw{\opm{\varphi}}A_1$ is an \cat{A}-epimorphism. With the help of this fact, we show that the morphism $(f,\varphi)$ is initial.

 \par
 Take a $\set\times\loc{A}$-morphism $|(X_3,\kappa_3,A_3)|\arw{(g,\psi)}|(X_1,\kappa_1,A_1)|$ with the property that $|(X_3,\kappa_3,A_3)|\arw{(g,\psi)}|(X_1,\kappa_1,A_1)|\arw{|(f,\varphi)|}|(X_2,\kappa_2,A_2)|$ is an \topsysto{L}-morphism. Consider the diagram
 $$
   \xymatrix{A_2\ar@/^1.5pc/"1,5"^-{\opm{(\varphi\circ\psi)}} \ar[d]_-{\opm{\kappa}_2}\ar[rr]_-{\opm{\varphi}} &&  A_1
             \ar[d]^-{\opm{\kappa}_1} \ar[rr]_-{\opm{\psi}} && A_3 \ar[d]^-{\opm{\kappa}_3}\\
             L^{X_2} \ar@/_1.5pc/"2,5"_-{\opm{(\mcal{P}_L(f\circ g))}} \ar[rr]^-{\opm{(\mcal{P}_Lf)}} && L^{X_1} \ar[rr]^-{\opm{(\mcal{P}_Lg)}} && L^{X_3}.}
 $$
 Commutativity of the outer square provides $\opm{\kappa}_3\circ\opm{\psi}\circ\opm{\varphi}=\opm{(\mcal{P}_Lg)}\circ\opm{\kappa}_1\circ\opm{\varphi}$ and, therefore,
 $\opm{\kappa}_3\circ\opm{\psi}=\opm{(\mcal{P}_Lg)}\circ\opm{\kappa}_1$ (since $\opm{\varphi}$ is an \cat{A}-epimorphism).
 \qed
\end{pf}

\begin{rem}
 \label{rem:1}
 Replacing Proposition~\ref{prop:11} by Proposition~\ref{prop:7} in the proof of Proposition~\ref{prop:12}, one arrives at the conclusion that in the category \topsyst{L}, embeddings are exactly the monomorphisms.
 \sqed
\end{rem}

\par
We recall that S.~Vickers~\cite[Definition~6.1.2]{Vickers1989} calls a topological system morphism $(X_1,\kappa_1,A_1)\arw{(f,\varphi)}(X_2,\kappa_2,A_2)$ an \emph{embedding} provided that $f$ (which is a map) is injective and \opm{\varphi} (which is a frame homomorphism) is surjective. It is well-known, however, that epimorphisms are surjective neither in the varieties \cat{Quant}, \cat{UQuant}~\cite{Kruml2008} nor in \cat{Frm}~\cite{Picado2012}. We define thus \mcal{M} to be the class of \topsysto{L}-monomorphisms $(f,\varphi)$ with \opm{\varphi} surjective, and characterize \mcal{M}-injective objects in the category \topsysto{L}. We notice that (by Proposition~\ref{prop:12}) \mcal{M} is a subclass of the class of embeddings in \topsysto{L}. In particular, if epimorphisms in \cat{A} are onto (e.g., in case of varieties $\semilatc{\bigvee}$ or
\cat{CBAlg}~\cite{Banaschewski2010}), then the two classes coincide.

\par
For convenience of the reader, we begin with recalling the definition of \mcal{M}-injective object in a category~\cite{Adamek2009}.

\begin{defn}
 \label{defn:6.1}
 Let \mcal{M} be a class of morphisms in a category \cat{C}. A \cat{C}-object $C$ is called \emph{\mcal{M}-injective} provided that for every \cat{C}-morphism $A\arw{m}B$ in \mcal{M} and every \cat{C}-morphism $A\arw{f}C$, there exists a \cat{C}-morphism $B\arw{g}C$, which makes the triangle
 $$
   \xymatrix{A \ar[r]^-{m} \ar[rd]_-{f} & B \ar[d]^-{g}\\
                                                  & C}
 $$
 commute.
 \sqed
\end{defn}

\begin{prop}
 \label{prop:13}
 $\mathsf{S}$ is an \mcal{M}-injective object in \topsysto{L}.
\end{prop}
\begin{pf}
 Given some $(X_1,\kappa_1,A_1)\arw{(f,\varphi)}(X_2,\kappa_2,A_2)\in\mcal{M}$ and an \topsysto{L}-morphism $(X_1,\kappa_1,A_1)\arw{(g,\psi)}\msf{S}$, we have to show the existence of a (necessarily unique) \topsysto{L}-morphism $(X_2,\kappa_2,A_2)\arw{(h,\theta)}\msf{S}$, making the following diagram commute
 $$
   \xymatrix{(X_1,\kappa_1,A_1) \ar[rr]^-{(f,\varphi)} \ar"2,3"_-{(g,\psi)} && (X_2,\kappa_2,A_2) \ar@{.>}[d]^-{(h,\theta)}\\
             & & \msf{S}.}
 $$

 \par
 Define $a=|\opm{\psi}|\circ\eta(\ast)$ and notice that since \opm{\varphi} is onto (by our assumption), there exists $b\in A_2$ such that $\opm{\varphi}(b)=a$. By Proposition~\ref{prop:2}, $(f_b,\varphi_b)\in\topsysto{L}((X_2,\kappa_2,A_2),\msf{S})$. By Proposition~\ref{prop:12}, it will be enough to verify that $(f_b,\varphi_b)\circ(f,\varphi)=(g,\psi)$.

 \par
 To show that $\varphi_b\circ\varphi=\psi$, we recall first (from Proposition~\ref{prop:2}) that $|\opm{\varphi}_b|\circ\eta=h_b^{A_2}$. Thus, $|\opm{\varphi}\circ\opm{\varphi}_b|\circ\eta(\ast)=\opm{\varphi}(b)=a=|\opm{\psi}|\circ\eta(\ast)$, which implies then $\opm{\varphi}\circ\opm{\varphi}_b=\opm{\psi}$.

 \par
 To show that $f_b\circ f=g$, we consider the following diagram
 $$
   \xymatrix{A_2\ar[d]_-{\opm{\kappa}_2}\ar[rr]^-{\opm{\varphi}} && A_1 \ar[d]^-{\opm{\kappa}_1}&& \ar[ll]_-{\opm{\psi}}\msf{S}
             \ar[d]^-{\opm{\kappa}_{S}} \ar@/_1.5pc/"1,1"_-{\opm{\varphi}_b}\\
             L^{X_2} \ar[rr]^-{\opm{(\mcal{P}_Lf)}} && L^{X_1} && \ar[ll]_-{\opm{(\mcal{P}_Lg)}} L^{|L|}, \ar@/^1.5pc/"2,1"^-{\opm{(\mcal{P}_Lf_b)}}}
 $$
 all the parts of which (except the lowest one) commute. Then, for every $x\in X_1$, $f_b\circ f(x)=(\opm{\kappa}_2(b))(f(x))=(\opm{(\mcal{P}_Lf)}\circ\opm{\kappa}_2(b))(x)=(\opm{\kappa}_1\circ\opm{\varphi}(b))(x)=
 (\opm{\kappa}_1(a))(x)=(\opm{\kappa}_1\circ\opm{\psi}(\eta(\ast)))(x)=(\opm{(\mcal{P}_Lg)}\circ\opm{\kappa}_S (\eta(\ast)))(x)=(\opm{\kappa}_S(\eta(\ast)))(g(x))=(|\pi_{g(x)}\circ\opm{\kappa}_S|\circ\eta)(\ast)=h_{g(x)}^L(\ast)=g(x)$.
 \qed
\end{pf}

\par
The following provides the main result of this subsection on the characterization of \mcal{M}-injective objects.

\begin{thm}
 \label{thm:3}
 \mcal{M}-injective objects in the category \topsysto{L} are precisely the retracts of powers of $\mathsf{S}$.
\end{thm}
\begin{pf}
 \hfill\par
 ``$\Longrightarrow$": We notice that given an \mcal{M}-injective $T_0$ affine system $(X,\kappa,A)$, there exists (by Theorem~\ref{thm:2}) an embedding $(X,\kappa,A)\arw{(f,\varphi)}\msf{S}^\mcal{F}$, in which $\mcal{F}=\topsysto{L}((X,\kappa,A),\msf{S})$.

 \par
 We first verify that $(f,\varphi)\in\mcal{M}$, i.e., $\opm{\varphi}$ is surjective. By Proposition~\ref{prop:2}, the elements of \mcal{F} are in one-to-one correspondence with the elements of $A$, and thus $\msf{S}^{\mcal{F}}=(|L|^{|A|},\theta,\copwr{|A|}{S})$. Moreover, the \cat{A}-homomorphism $\copwr{|A|}{S}\arw{\opm{\varphi}}A$ is given by the diagram (which is the universal property of coproducts)
 $$
   \xymatrix{S \ar[rd]_-{\opm{\varphi}_a} \ar[r]^-{\mu_a} & \copwr{|A|}{S} \ar@{.>}[d]^-{\opm{\varphi}}\\
             & A.}
 $$
 Given $a\in A$, $a=h_a^A(\ast)=|\opm{\varphi}_a|\circ\eta(\ast)=
 \opm{\varphi}_a(\eta(\ast))=\opm{\varphi}\circ\mu_a(\eta(\ast))=
 \opm{\varphi}(\mu_a(\eta(\ast)))$, i.e., $\opm{\varphi}$ is onto.

 \par
 Since $(X,\kappa,A)$ is \mcal{M}-injective, there exists an \topsysto{L}-morphism $\msf{S}^\mcal{F}\arw{(g,\psi)}(X,\kappa,A)$, which makes the following diagram commute
 $$
   \xymatrix{(X,\kappa,A) \ar[rd]_-{1_{(X,\kappa,A)}} \ar[r]^-{(f,\varphi)} & \msf{S}^{\mcal{F}} \ar@{.>}[d]^-{(g,\psi)}\\
             &(X,\kappa,A).}
 $$
 It follows then that $\msf{S}^\mcal{F}\arw{(g,\psi)}(X,\kappa,A)$ is a retraction.

 \par
 ``$\Longleftarrow$": We notice that retracts of \mcal{M}-injective objects are \mcal{M}-injective~\cite[Proposition~9.5]{Adamek2009}, and products of \mcal{M}-injective objects are \mcal{M}-injective~\cite[Proposition~10.40]{Adamek2009}.
 \qed
\end{pf}

\par
As a particular case of Theorem~\ref{thm:3}, we obtain the result of~\cite{Noor2016} on the characterization of \mcal{M}-injective $T_0$ topological systems as retracts of powers of the Sierpinski topological system. Moreover, an attentive reader will probably notice that Theorem~\ref{thm:3} can be also shown without the assumptions on the variety \cat{A} and its algebra $L$, mentioned at the very beginning of this subsection.


\subsection{Sober affine systems}


\par
In this subsection, we will characterize sober affine systems with the help of the Sierpinski affine system. We begin with an affine modification of the concept of localic system of~\cite{Vickers1989}. Given an \cat{A}-algebra $A$, from now on, we will employ the notation $Pt_L(A)=\cat{A}(A,L)$. The elements of $Pt_L(A)$ will be denoted $p$ (``$p$" being an abbreviation for ``point" as motivated by, e.g.,~\cite{Johnstone1982}).

\begin{prop}
 \label{prop:14}
 Given an affine system $(X,\kappa,A)$, there exists a map $X\arw{\ell}Pt_L(A)$, $\ell(x)=(\opm{\kappa}(-))(x)$.
\end{prop}
\begin{pf}
 Given $x\in X$, we have to verify that $\ell(x)\in Pt_L(A)$. Given $\lambda\in\Lambda$ (we employ here the notations of Definition~\ref{defn:1}), it follows that $(\ell(x))(\omega_{\lambda}^A(\indf{a}{i}))=(\opm{\kappa}(\omega_{\lambda}^A(\indf{a}{i})))(x)=
 (\omega_{\lambda}^{L^X}(\indf{\opm{\kappa}(a_i)}{}))(x)=\omega_{\lambda}^L(\indf{(\opm{\kappa}(a_i))(x)}{})=
 \omega_{\lambda}^L(\indf{(\ell(x))(a_i)}{})$.
 \qed
\end{pf}

\begin{defn}
 \label{defn:7}
 An affine system $(X,\kappa,A)$ is \emph{sober} provided that the map $X\arw{\ell}Pt_L(A)$ is bijective.
 \sqed
\end{defn}

\par
For convenience of the reader, we notice that in case of the embedding \incl{\cat{Top}}{E}{\cat{TopSys}} of Theorem~\ref{thm:1}, sober systems of the form $E(X,\tau)$ are precisely sober topological spaces~\cite{Johnstone1982} (thus the term ``sober"). We also emphasize that the definition of $T_0$ affine systems (Definition~\ref{defn:6}) implies immediately the next result.

\begin{prop}
 \label{prop:15}
 An affine system $(X,\kappa,A)$ is $T_0$ iff the map $X\arw{\ell}Pt_L(A)$ is injective. In particular, every sober affine system is $T_0$.
\end{prop}

\par
We now provide an (probably, already expected by the reader) example of sober affine systems.

\begin{prop}
 \label{prop:16}
 The Sierpinski affine system \msf{S} is sober.
\end{prop}
\begin{pf}
 We have to show that the map $|L|\arw{\ell}Pt_L(S)$ is bijective. By Proposition~\ref{prop:4}, \msf{S} is $T_0$, which implies injectivity of $\ell$ (Proposition~\ref{prop:15}). To show surjectivity, we fix $p\in Pt_L(S)$ and define $a=|p|\circ\eta(\ast)$. For every $b\in S$, $(\ell(a))(b)=(\opm{\kappa}_S(b))(a)=
 \pi_a\circ\opm{\kappa}_S(b)=\ovr{h^L_a}(b)\overset{(\dagger)}{=}p(b)$, where $(\dagger)$ relies on the fact that $|p|\circ\eta(\ast)=a=|\ovr{h^L_a}|\circ\eta(\ast)$ implies $p=\ovr{h^L_a}$. It follows then that $\ell(a)=p$, namely, $\ell$ is surjective.
 \qed
\end{pf}

\par
To arrive at the main theorem of this subsection, we continue with several preliminary results.

\begin{prop}
 \label{prop:17}
 Every sober affine system can be embedded into some power of \msf{S}.
\end{prop}
\begin{pf}
 Follows from Proposition~\ref{prop:15} and Theorem~\ref{thm:2}.
 \qed
\end{pf}

\par
We are going to show that the embedding of Proposition~\ref{prop:17} has an additional property, which characterizes sober affine systems among other affine systems.

\begin{prop}
 \label{prop:18}
 The product of a set-indexed family of sober affine systems is a sober affine system.
\end{prop}
\begin{pf}
 Given a set-indexed family $\{(X_i,\kappa_i,A_i)\mid i\in I\}$ of sober affine systems, we show that the product $(\prod_{i\in I}X_i,\kappa,\coprod_{i\in I}A_i)$ is sober. Consider the map $\prod_{i\in I}X_i\arw{\ell}Pt_L(\coprod_{i\in I}A_i)$. Since products of $T_0$ systems are $T_0$ (Proposition~\ref{prop:9}), $\ell$ is injective (Proposition~\ref{prop:15}). To show that the map is also surjective, we take some $p\in Pt_L(\coprod_{i\in I}A_i)$. For every $i\in I$, consider the commutative diagram
 $$
   \xymatrix{A_i \ar[rd]_-{p_i}\ar[r]^-{\mu_i}&\coprod_{i\in I}A_i \ar[d]^-{p}\\
             & L.}
 $$
 For every $i\in I$, since $(X_i,\kappa_i,A_i)$ is sober, there exists $x_i\in X_i$ such that $\ell_i(x_i)=p_i$. Thus, given $j\in I$ and $a\in A_j$, $(\ell(\indfn{x}{i}{i\in I}))\circ\mu_j(a)=(\opm{\kappa}(\mu_j(a)))(\indfn{x}{i}{i\in I})=
 (\opm{\kappa}\circ\mu_j(a))(\indfn{x}{i}{i\in I})=(\opm{(\mcal{P}_L\pi_j)}\circ\opm{\kappa}_j(a))(\indfn{x}{i}{i\in I})=
 (\opm{\kappa}_j(a))(\pi_j(\indfn{x}{i}{i\in I}))=(\opm{\kappa}_j(a))(x_j)=(\ell_j(x_j))(a)=p_j(a)$. It follows then that $\ell(\indfn{x}{i}{i\in I})=p$.
 \qed
\end{pf}

\par
We provide now an equivalent description of the condition on affine morphisms from Definition~\ref{defn:4}.

\begin{prop}
 \label{prop:19}
 Given two affine systems $(X_1,\kappa_1,A_1)$ and $(X_2,\kappa_2,A_2)$, and a $\set\times\loc{A}$-morphism $|(X_1,\kappa_1,A_1)| \arw{(f,\varphi)}|(X_2,\kappa_2,A_2)|$, $(f,\varphi)$ is an affine morphism iff the following diagram commutes (its lower arrow is actually a restriction of the map $L^{|A_1|}\arw{\opm{(\mcal{P}_L|\opm{\varphi}|)}}L^{|A_2|}$, since $Pt_L(A_i)$ is a subset of $L^{|A_i|}$)
 \begin{equation}
  \tag{D}
  \label{diag:1}
  \xymatrix{X_1 \ar[d]_{\ell_1} \ar[rr]^{f} & & X_2 \ar[d]^{\ell_2}\\
            Pt_L(A_1) \ar[rr]_{|\opm{(\mcal{P}_L|\opm{\varphi}|)}|} && Pt_L(A_2).}
 \end{equation}
\end{prop}
\begin{pf}
 For the necessity, we notice that given $x_1\in X_1$ and $a_2\in A_2$, $(|\opm{(\mcal{P}_L|\opm{\varphi}|)}|\circ\ell_1(x_1))(a_2)=(\ell_1(x_1)\circ\opm{\varphi})(a_2)=(\ell_1(x_1))(\opm{\varphi}(a_2))=
 (\opm{\kappa}_1(\opm{\varphi}(a_2)))(x_1)=(\opm{\kappa}_1\circ\opm{\varphi}(a_2))(x_1)=
 (\opm{(\mcal{P}_Lf)}\circ\opm{\kappa}_2(a_2))(x_1)=(\opm{\kappa}_2(a_2))(f(x_1))=(\ell_2\circ f(x_1))(a_2)$. For the sufficiency, fixing again $x_1\in X_1$ and $a_2\in A_2$, we notice that $(\opm{(\mcal{P}_Lf)}\circ\opm{\kappa}_2(a_2))(x_1)=(\opm{\kappa}_2(a_2))(f(x_1))=(\ell_2\circ f(x_1))(a_2)=(\opm{|(\mcal{P}_L|\opm{\varphi}|)}|\circ\ell_1(x_1))(a_2)=(\ell_1(x_1))(\opm{\varphi}(a_2))=
 (\opm{\kappa}_1(\opm{\varphi}(a_2)))(x_1)=(\opm{\kappa}_1\circ\opm{\varphi}(a_2))(x_1)$.
 \qed
\end{pf}

\par
In view of Proposition~\ref{prop:19}, we introduce the following notion.

\begin{defn}
 \label{defn:8}
 An affine monomorphism $(X_1,\kappa_1,A_1)\arw{(f,\varphi)}(X_2,\kappa_2,A_2)$ is \emph{sober} provided that (\ref{diag:1}) is a weak pullback~\cite{Hofmann2014} (namely, the canonical map $X_1\arw{}Pt_L(A_1)\times_{Pt_L(A_2)}X_2$ is surjective).
 \sqed
\end{defn}

\par
For convenience of the reader, we provide a brief comment on Definition~\ref{defn:8} w.r.t. affine spaces. We recall first from~\cite{Solovyov2015e} a convenient property of the category \topt{L}.

\begin{thm}
 \label{thm:4}
 The concrete category $(\topt{L}$, ${|-|})$ is topological over \set.
\end{thm}
\begin{pf}
 Given a $|-|$-structured source $\mcal{L}=(X\arw{f_i}|(X_i,\tau_i)|)_{i\in I}$, the initial structure on the set $X$ w.r.t. \mcal{L} is given by the subalgebra of $L^X$, which is generated by the union $\bigcup_{i\in I}\opm{(\mcal{P}_Lf_i)}(\tau_i)$.
 \qed
\end{pf}

\par
The following corollary of Theorem~\ref{thm:4} employs the embedding functor of Theorem~\ref{thm:1}.

\begin{cor}
 \label{cor:1}
 Let $(X_1,\tau_1)\arw{f}(X_2,\tau_2)$ be an \topt{L}-morphism. If $f$ is an embedding, then $Ef$ is an embedding (the same as monomorphism by Remark~\ref{rem:1}) in \topsyst{L}. The converse though does not hold.
\end{cor}
\begin{pf}
 For the first statement, we notice that since $(X_1,\tau_1)\arw{f}(X_2,\tau_2)$ is an embedding, $f$ is injective. By Theorem~\ref{thm:1},  $E((X_1,\tau_1)\arw{f}(X_2,\tau_2))=
 (X_1,\opm{e}_{\tau_1},\tau_1)\arw{(f,\varphi)}(X_2,\opm{e}_{\tau_2},\tau_2)$, with $e_{\tau_i}$ the inclusion $\tau_i\hookrightarrow L^{X_i}$ and \opm{\varphi} the restriction $\tau_2\arw{\opm{(\mcal{P}_Lf)}|_{\tau_2}^{\tau_1}}\tau_1$. By the proof of Theorem~\ref{thm:4}, \opm{\varphi} is an epimorphism in \cat{A}, and, therefore, $Ef$ is an \topsyst{L}-monomorphism, i.e., an embedding in \topsyst{L}.

 \par
 For the second statement, we show the following two simple counterexamples. First, \cite[2.9.3]{Barr1990} states that the inclusion $\mbb{N}\overset{e}{\hookrightarrow}\mbb{Z}$ of natural numbers into integers is an epimorphism in the variety \cat{Mon} of monoids, where $e(\mbb{N})=\mbb{N}$ is not the whole \mbb{Z}. If $\cat{A}=\cat{Mon}$ and $L=\mbb{Z}$, then for $X_1=X_2=\msf{1}$, $\tau_1=L^{\msf{1}}$, $\tau_2=\{\alpha\in L^{\msf{1}}\mid\alpha(\ast)\in\mbb{N}\}$, and $X_1\arw{!}X_2$ the unique possible map, $E!$ is an embedding in \topsyst{L}, but $(X_1,\tau_1)\arw{!}(X_2,\tau_2)$ is not an embedding in \topt{L}. Second, \cite[Example~6.1.1]{Picado2012} states that in the category \cat{Top} (thus, $\cat{A}=\cat{Frm}$ and $L=\msf{2}$), given a $T_1$-space $(X,\tau)$, the identity map $(X,\mcal{P}X)\arw{1_X}(X,\tau)$ provides a \cat{Frm}-epimorphism $\tau\hookrightarrow\mcal{P}X$, which is surjective iff $\tau=\mcal{P}X$. Therefore, $E1_X$ is an embedding in \cat{TopSys}, but $(X,\mcal{P}X)\arw{1_X}(X,\tau)$ is not an embedding in \cat{Top} provided that $\tau\not=\mcal{P}X$.
 \qed
\end{pf}

\par
As follows from Corollary~\ref{cor:1}, the concept of embedding in the category \topsyst{L} is ``strictly weaker" than the concept of embedding in the category \topt{L}. Further, we notice that in the classical setting of the categories \cat{Top} and \cat{TopSys}, given a continuous map $(X_1,\tau_1)\arw{f}(X_2,\tau_2)$ such that $Ef$ is a sober monomorphism (the latter being the same as embedding), we can only assume that, first, $X_1\seq X_2$, i.e., $X_1\arw{f}X_2$ is the inclusion $X_1\hookrightarrow X_2$, and, second,  $\tau_2\arw{(\opm{\mcal{P}f)}|_{\tau_2}^{\tau_1}}\tau_1$, $\opm{(\mcal{P}f)}(U)=U\bigcap X_1$ is a \cat{Frm}-epimorphism (which need not be surjective by the proof of Corollary~\ref{cor:1}). The condition of Definition~\ref{defn:8} then states (cf., e.g.,~\cite{Picado2012,Picado2004}) that given a completely prime filter (cp-filter, for short) $\mfrak{F}_1\in Pt(\tau_1)$, if the cp-filter $\mfrak{F}_2=\{U\in\tau_2\mid U\bigcap X_1\in\mfrak{F}_1\}$ is precisely the neighborhood filter $\mfrak{U}_2(y)=\{U\in\tau_2\mid y\in U\}$ for some $y\in X_2$, then $\mfrak{F}_1=\mfrak{U}_1(x)=\{U\in\tau_1\mid x\in U\}$ for some $x\in X_1$. Moreover, in case of a $T_0$-space $(X_2,\tau_2)$, one gets that $x=y$, namely, $y\in X_1$. In case of a sober space $(X_2,\tau_2)$, $X_1$ must contain every $y\in X_2$ with the property that $\mfrak{U}_2(y)=\{U\in\tau_2\mid U\bigcap X_1\in\mfrak{F}_1\}$ for some $\mfrak{F}_1\in Pt(\tau_1)$.

\begin{prop}
 \label{prop:20}
 Sober subobjects of sober affine systems are sober.
\end{prop}
\begin{pf}
 Given a sober \topsyst{L}-monomorphism $(X_1,\kappa_1,A_1)\arw{(f,\varphi)}(X_2,\kappa_2,A_2)$ such that $(X_2,\kappa_2,A_2)$ is sober, we have to show that $(X_1,\kappa_1,A_1)$ is also sober. Since every sober affine system is $T_0$, Propositions~\ref{prop:10}, \ref{prop:15} imply that the map $X_1\arw{\ell_1}Pt_L(A_1)$ is injective. To show that the map is also surjective, we fix $p_1\in Pt_L(A_1)$. Then $|\opm{(\mcal{P}_L|\opm{\varphi}|)}|(p_1)=p_2\in Pt_L(A_2)$, and, therefore, there exists $x_2\in X_2$ such that $\ell_2(x_2)=p_2$. Since Diagram~(\ref{diag:1}) is a weak pullback, there exists $x_1\in X_1$ such that $\ell_1(x_1)=p_1$.
 \qed
\end{pf}

\par
We provide now the main result of this subsection on the characterization of sober affine systems. Our employed term ``soberly embedded" means that the embedding in question is sober (cf. Definition~\ref{defn:8}).

\begin{thm}
 \label{thm:5}
 An affine system $(X,\kappa,A)$ is sober iff it can be soberly embedded into some power of $\mathsf{S}$.
\end{thm}
\begin{pf}
 \hfill\par
 ``$\Longleftarrow$": Follows from Propositions~\ref{prop:16}, \ref{prop:18}, \ref{prop:20}.

 \par
 ``$\Longrightarrow$": By Theorem~\ref{thm:2}, there exists an embedding $(X,\kappa,A)\arw{(f,\varphi)}\msf{S}^\mcal{F}$ with $\mcal{F}=\topsysto{L}((X,\kappa,A),\msf{S})$. Moreover, we have already mentioned in the proof of Theorem~\ref{thm:3} that the elements of \mcal{F} are in one-to-one correspondence with the elements of $A$, and thus $\msf{S}^{\mcal{F}}=(|L|^{|A|},\theta,\copwr{|A|}{S})$. The only thing left to verify is that the following commutative (by Proposition~\ref{prop:19}) diagram is a weak pullback
 $$
   \xymatrix{X \ar[d]_{\ell} \ar[rr]^{f} & & |L|^{|A|} \ar[d]^{\ell_S}\\
             Pt_L(A) \ar[rr]_{|\opm{(\mcal{P}_L|\opm{\varphi}|)}|} && Pt_L(\copwr{|A|}{S}).}
 $$

 \par
 Fix $p\in Pt_L(A)$ and $\indfn{b}{a}{A}\in|L|^{|A|}$ such that $|\opm{(\mcal{P}_L|\opm{\varphi}|)}|(p)=\ell_S(\indfn{b}{a}{A})$. Since $(X,\kappa,A)$ is sober, there exists $x\in X$ such that $\ell(x)=p$. We show that $f(x)=\indfn{b}{a}{A}$, which will finish the proof. Indeed, $\ell_S\circ f(x)=(|\opm{(\mcal{P}_L|\opm{\varphi}|)}|\circ\ell)(x)=
 |\opm{(\mcal{P}_L|\opm{\varphi}|)}|(p)=\ell_S(\indfn{b}{a}{A})$. Since $\ell_S$ is injective, $f(x)=\indfn{b}{a}{A}$.
 \qed
\end{pf}

\par
We notice that Theorem~\ref{thm:5} goes beyond the results of~\cite{Noor2016} and is motivated by the results in~\cite{Nel1972}. More precisely, L.~D.~Nel and R.~G.~Wilson considered in~\cite{Nel1972} the so-called \emph{front topology} on a topological space $(X,\tau)$ by specifying the \emph{front-closure} operator \emph{fcl} as follows: given $x\in X$ and $Y\seq X$, $x\in fcl(Y)$ iff $U\bigcap cl(\{x\})\bigcap Y\not=\varnothing$ for every $U\in\mfrak{U}(x)$ (one notices that for every non-discrete $T_0$-space, the front topology is strictly larger than the original topology). Theorem~3.4 in~\cite{Nel1972} states then that a $T_0$-space is sober iff it is homeomorphic to a front-closed subspace of some power of the Sierpinski space. Our remark, following Corollary~\ref{cor:1}, provides a topological system analogue of the above-mentioned closure operator \emph{fcl}.


\section{\label{sec:1}Sierpinski space versus Sierpinski system}


\par
In this section, we compare the Sierpinski affine space and the Sierpinski affine system. We notice first that there already exists a lattice-valued analogue of the Sierpinski space~\cite{Srivastava1984,Srivastava1986}. Moreover, its affine version has already been studied in, e.g.,~\cite{Singh2013,Solovyov2008b}, which motivates our next definition.

\begin{defn}
 \label{defn:9}
 \emph{Sierpinski ($L$-)affine space} is the pair $\mcal{S}=(|L|,\subb{1_L})$, where $\subb{1_L}$ stands for the subalgebra of $L^{|L|}$, which is generated by the identity map $1_L$.
 \sqed
\end{defn}

\par
According to~\cite[Theorem~3.2]{Singh2013} (cf. also~\cite[Theorem~58]{Solovyov2008b}), $\mcal{S}$ is a Sierpinski object in the category \topt{L}, which we will call the Sierpinski affine space as in Definition~\ref{defn:9}.

\par
For convenience of the reader, we show an example of the Sierpinski affine space, kindly suggested by one of the referees. We recall that a quantale $(L,\otimes)$ is \emph{idempotent} provided that $a\otimes a=a$ for every $a\in L$.

\begin{exmp}
 \label{exmp:5}
 Let $\cat{A}=\cat{UQuant}$ and let $L$ be an integral quantale (i.e., $\unit_L=\top_L$). The respective Sierpinski affine space then has the form $\mcal{S}=(|L|,\subb{1_L}=\{\underline{\bot_L},\underline{\top_L}\}\cup\{(1_L)^n\mid n\in\mbb{N}\})$, where \mbb{N} is the set of natural numbers $\{1,2,3,\ldots\}$ and $(1_L)^n=\underbrace{1_L\otimes\ldots\otimes 1_L}_{\text{$n$-times}}$. Moreover, if $L$ is idempotent, then $\subb{1_L}=\{\underline{\bot_L},1_L,\underline{\top_L}\}$. In particular, if $\cat{A}=\cat{Frm}$ and $L=\msf{2}$, one gets the classical Sierpinski space from the introductory section.
 \sqed
\end{exmp}

\par
One can also show the following result~\cite{Solovyov2008b} (which is exactly item~(1) of the introductory section).

\begin{defn}
 \label{defn:10}
 An affine space $(X,\tau)$ is $T_0$ provided that for every $x,y\in X$, $\alpha(x)=\alpha(y)$ for very $\alpha\in\tau$ implies $x=y$.
 \sqed
\end{defn}

\begin{prop}
 \label{prop:21}
 An affine space $(X,\tau)$ is $T_0$ iff it can be embedded into some power of \mcal{S}.
\end{prop}

\par
It is not our purpose to study the properties of the Sierpinski affine space (which could be the topic of our forthcoming papers), but rather to compare it with the Sierpinski affine system. More precisely, considering Theorem~\ref{thm:1}, one can ask the question  whether $E\mcal{S}$ is ``comparable" (e.g., isomorphic) to \msf{S}. In the following, we try to give a partial answer to this question.

\par
We notice first that there clearly exists an \cat{A}-homomorphism $S\arw{\ovr{h_{1_L}^{\subb{1_L}}}}\subb{1_L}$ (for the sake of brevity, denoted $\vartheta$), which is given by the diagram
$$
  \xymatrix{\msf{1} \ar[rd]_-{h_{1_L}^{\subb{1_L}}}\ar[r]^-{\eta}& |S| \ar@{.>}[d]^-{|\vartheta|}\\
             & |\subb{1_L}|.}
$$
Thus, there exists a $\set\times\loc{A}$-morphism $|E\mcal{S}=(|L|,\opm{e}_{\subb{1_L}},\subb{1_L})|\arw{(1_L,\opm{\vartheta})}|\msf{S}=(|L|,\kappa_S,S)|$.

\begin{prop}
 \label{prop:22}
 The $\set\times\loc{A}$-morphism $|E\mcal{S}|\arw{(1_L,\opm{\vartheta})}|\msf{S}|$ is an \topsyst{L}-morphism.
\end{prop}
\begin{pf}
 We consider the following diagram
 $$
   \xymatrix{\msf{1}\ar@/_2.3pc/"3,4"_-{h_a^L}\ar@/^1.5pc/"1,4"^-{h_{1_L}^{\subb{1_L}}} \ar[r]^-{\eta}&|S|
             \ar[d]_-{|\opm{\kappa}_S|}\ar[rr]^-{|\vartheta|} && |\subb{1_L}|
             \ar[d]^-{|\opm{e}_{\subb{1_L}}|}&\\
             &|L^{|L|}| \ar[rr]_-{|1_{L^{|L|}}|} && |L^{|L|}| \ar[d]^-{|\pi_a|}\\
             &&&|L|}
 $$
 and notice that $|\pi_a\circ\opm{e}_{\subb{1_L}}\circ\vartheta|\circ\eta(\ast)=|\pi_a\circ\opm{e}_{\subb{1_L}}\circ h_{1_L}^{\subb{1_L}}|(\ast)=|\pi_a|(1_L)=|1_L|(a)=a=h_a^L(a)=|\pi_a\circ 1_{L^{|L|}}\circ\opm{\kappa}_S|\circ\eta(\ast)$ for every $a\in L$. Since the products in the variety \cat{A} are concrete (i.e., are preserved by the forgetful functor), we get
 $|\opm{e}_{\subb{1_L}}\circ\vartheta|\circ\eta(\ast)=|1_{L^{|L|}}\circ\opm{\kappa}_S|\circ\eta(\ast)$ and then $\opm{e}_{\subb{1_L}}\circ\vartheta=1_{L^{|L|}}\circ\opm{\kappa}_S$.
 \qed
\end{pf}

\par
Proposition~\ref{prop:22} implies immediately the following corollary.

\begin{cor}
 \label{cor:2}
 The affine systems $E\mcal{S}$ and \msf{S} are isomorphic provided that $S\arw{\vartheta}\subb{1_L}$ is an \cat{A}-isomorphism.
\end{cor}

\par
In the following, we check, whether Corollary~\ref{cor:2} is true in case of the varieties of Example~\ref{exmp:1}. Given two sets $X$, $Y$ and an element $y\in Y$, we will denote by $\underline{y}$ the constant map $X\arw{}Y$ with value $y$.
\begin{description}
 \item[\semilatc{\bigvee}:] The free $\bigvee$-semilattice $S$ over a singleton is the two-element chain
       $\{\bot_S,\top_S\}$ (with $\eta(\ast)=\top_S$). If the $\bigvee$-semilattice $L$ has more than one element, then $\subb{1_L}=\{\underline{\bot_L},1_L\}$. Thus, the map $S\arw{\vartheta}\subb{1_L}$ is an isomorphism of $\bigvee$-semilattices. As a consequence, Corollary~\ref{cor:2} holds in case of $\bigvee$-semilattices (recall that starting from Subsection~\ref{subsec:1}, we require $L$ to have at least two elements).
 \item[\cat{(U)Quant}:] The free (resp. unital) quantale over a singleton is the powerset of the set $\mathbb{N}$ (resp.
      $\mathbb{N}\bigcup\{0\}$) of natural numbers~\cite{Rosenthal1990} and thus, is uncountably infinite. Following Example~\ref{exmp:5}, on the one hand, if $L$ is an integral idempotent quantale (e.g., a frame), then $\subb{1_L}$ has at most three elements, i.e., the map $S\arw{\vartheta}\subb{1_L}$ is \emph{not} an isomorphism of (resp. unital) quantales. On the other hand, if $L$ is an integral non-idempotent quantale (e.g., the real unit interval $[0,1]$ with the {\L}ukasiewicz t-norm~\cite{Hohle2011}), then $\subb{1_L}$ is at most countably infinite, i.e., the map $S\arw{\vartheta}\subb{1_L}$ is \emph{not} an isomorphism of (resp. unital) quantales. As a consequence, Corollary~\ref{cor:2} is false for (resp. unital) quantales. One can also see the latter statement as follows. Since the map $S\arw{\opm{\kappa}_S}L^{|L|}$ of Example~\ref{exmp:4}, considered as a map onto its range, coincides with the map $S\arw{\vartheta}\subb{1_L}$ (which follows from the universal property of free objects), and $\opm{\kappa}_S$ is not injective (cf. Example~\ref{exmp:4}), then so is $\vartheta$, i.e., $\vartheta$ is \emph{not} an isomorphism.
 \item[\cat{Frm}:] We have already mentioned that the free frame $S$ over a singleton is given by the three-element chain
      $\{\bot_S,c,\top_S\}$. If $L$ has more than one element (in particular, if $L=\msf{2}$), then $\subb{1_L}=\{\underline{\bot_L},1_L,\underline{\top_L}\}$. Thus, the map $S\arw{\vartheta}\subb{1_L}$ is an isomorphism of frames. As a consequence, Corollary~\ref{cor:2} holds in case of frames. In particular, the setting of $L=\msf{2}$ provides an isomorphism between $E\mcal{S}$ (the image under $E$ of the classical Sierpinski space) and $\msf{S}$ (Sierpinski topological system of~\cite{Noor2016}).
 \item[\cat{CBAlg}:] The free complete Boolean algebra $S$ over a singleton has the form
      $$
        
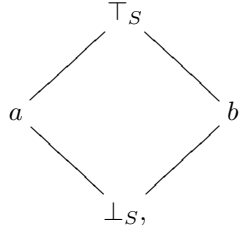
\xymatrix{&\ar@{-}[ld]^-{}\top_S \ar@{-}[rd]^-{}&\\
                  a\ar@{-}[rd]_-{}&&b \ar@{-}[ld]_-{}\\
                  &\bot_S,&}
      $$
      where $a^{\ast}=b$. If $L$ has more than one element (in particular, if $L=\msf{2}$), then $\subb{1_L}=\{\underline{\bot_L},1_L,1^{\ast}_L,\underline{\top_L}\}$ (where $1^{\ast}_L(a)=a^{\ast}$ for every $a\in L$). Thus, the map $S\arw{\vartheta}\subb{1_L}$ is an isomorphism of complete Boolean algebras. As a consequence, Corollary~\ref{cor:2} holds in case of complete Boolean algebras.
\end{description}

\par
As can be seen from the above discussion, in case of almost all varieties of Example~\ref{exmp:1} (quantales make the only exception), the form of the Sierpinski affine system is already predetermined by the form of the Sierpinski affine space. Thus, it seems plausible that the properties of the latter could be translated into the language of affine systems, making them the properties of the former.


\section{Conclusion}


\par
This paper makes another step in our effort to bring the theory of lattice-valued topology under the setting of affine sets of~Y.~Diers~\cite{Diers1996,Diers1999,Diers2002}. In particular, we have considered an affine setting for topological systems of S.~Vickers~\cite{Vickers1989} and introduced an affine system analogue of the well-known Sierpinski space. Our study was motivated by the paper of R.~Noor and A.~K.~Srivastava~\cite{Noor2016}, who presented the Sierpinski topological system and studied its basis properties. With the help of our affine setting, we have extended their results from crisp case to many-valued case and included additional results. For example, we showed that in the affine setting of frames (which includes the standard notions of topological space and topological system), the Sierpinski affine system is isomorphic to the Sierpinski affine space (when represented as a system). Such a result could potentially allow us to translate the results on the Sierpinski space into the language of topological systems (which could be the topic of our forthcoming papers). We would like to conclude the discussion of this paper with several open problems, which seem of interest to us.

\par
We have considered the simplest possible case of fixed-basis affine systems (and spaces) over the category \set, i.e., the category \topsyst{L} (and \topt{L}). We see here two possible directions for further study.

\begin{prob}
 \label{prob:1}
 Extend our fixed-basis results to the variable-basis case of Proposition~\ref{prop:1}.
 \sqed
\end{prob}

\begin{prob}
 \label{prob:2}
 Extend our \set-valued results to the case of general functor $\cat{X}\arw{T}\loc{A}$.
 \sqed
\end{prob}

\par
In this paper, we have restricted ourselves to just five varieties of Example~\ref{exmp:1}. As follows from the discussion at the end of Section~\ref{sec:1}, only two of them (the varieties of (unital) quantales) do not make Corollary~\ref{cor:2} true. Since the validity of the corollary in question could make it easier to study the properties of the Sierpinski affine system (through the properties of the Sierpinski affine space), we get our last problem.

\begin{prob}
 \label{prob:3}
 What are the conditions on the variety \cat{A}, which ensure the validity of Corollary~\ref{cor:2} (namely, which make the \cat{A}-homomorphism $S\arw{\vartheta}\subb{1_L}$ an isomorphism)?
 \sqed
\end{prob}


\section*{Acknowledgements}


\par
The authors would like to express their sincere gratitude to the anonymous referees of this paper for their helpful remarks and suggestions, e.g., Examples~\ref{exmp:4}, \ref{exmp:5} are due to one of them.

\par
This is a preprint of the paper published in ``Fuzzy Sets and Systems". The final authenticated version of the paper is available online at: https://www.sciencedirect.com/science/article/pii/S0165011416302706.



\end{document}